\documentclass{amsart}

\begin{document}

\title[Hasse polynomial]{Lie Invariant Frobenius lifts on\\ linear algebraic groups}

\def \h{\hat{\ }}
\def \cO{\mathcal O}
\def \ra{\rightarrow}
\def \bZ{{\mathbb Z}}
\def \cP{\mathcal V}
\def \cH{{\mathcal H}}
\def \cB{{\mathcal B}}
\def \d{\delta}
\def \cC{{\mathcal C}}
\def \jor{\text{jor}}

\newtheorem{THM}{{\!}}[section]
\newtheorem{THMX}{{\!}}
\renewcommand{\theTHMX}{}
\newtheorem{theorem}{Theorem}[section]
\newtheorem{corollary}[theorem]{Corollary}
\newtheorem{lemma}[theorem]{Lemma}
\newtheorem{proposition}[theorem]{Proposition}
\newtheorem{thm}[theorem]{Theorem}
\theoremstyle{definition}
\newtheorem{definition}[theorem]{Definition}
\theoremstyle{remark}
\newtheorem{remark}[theorem]{Remark}
\newtheorem{example}[theorem]{Examples}
\numberwithin{equation}{section}
\newtheorem{conjecture}[theorem]{\bf Conjecture}
\address{Department of Mathematics and Statistics\\University of New Mexico \\
  Albuquerque, NM 87131, USA\\ 
   Department of Mathematics} 
\subjclass[2010]{Primary: 20G25, Secondary: 11F85}
\maketitle

\bigskip

\medskip
\centerline{\bf Alexandru Buium}
\bigskip

\begin{abstract} 
We show  that if $G$ is a linear algebraic group over a number field and if $G$ is not a torus  then for all but finitely many primes $p$ the $p$-adic completion of $G$ does not possess a  Frobenius lift that is ``Lie invariant mod $p$" (in the sense of \cite{alie1}). This is in contrast 
with the situation of elliptic curves studied in
\cite{alie1}. \end{abstract}

\bigskip
\bigskip
\bigskip

\section{Introduction}
Consider a  complete discrete valuation ring $R$ with maximal ideal generated by a prime $p$ and algebraically closed residue field $\overline{R}:=R/pR$. Then   Frobenius lifts
 on $p$-adic formal
schemes  over  $R$ can be viewed as arithmetic analogues of vector fields on complex algebraic varieties;  cf.  \cite{char, book, foundations}.  In view of this analogy one is naturally led to ask the following question: {\it what  Frobenius lifts on $p$-adic completions of  group schemes over $R$ correspond to (translation)  invariant vector fields on complex algebraic groups}? A possible answer to this question is provided by the  following concept. 
 Let $G$ be a smooth  group scheme over $R$. Let us say that a Frobenius lift $\phi$ on the $p$-adic completion $\widehat{G^0}$ of an open subset $G^0\subset G$ is   {\it Lie invariant mod $p$}  if 
 \begin{equation}
\label{inver}
\frac{\phi^*}{p}\omega\equiv \lambda \omega\ \ \ \text{mod}\ \ \ p,\end{equation}
where $\omega$ is a column whose entries form an $R$-basis for the left invariant $1$-forms on $G$ and $\lambda$ is a matrix with $R$-coefficients. This condition is independent of $\omega$
and can be viewed as an arithmetic analogue of the condition that the Lie derivative along a vector field on a complex algebraic group kills the left invariant forms. 
 A similar definition can be made mod $p^m$. (Other arithmetic analogues  of invariant vector fields were considered in \cite{foundations}; we will not discuss them here.)
 
 I view of \cite{alie1} 
 one is tempted to conjecture that any ordinary Abelian scheme $G$ has an open set $G^0$ whose $p$-adic completion $\widehat{G^0}$  possesses a 
    Frobenius lift that is Lie invariant mod $p$. 
    In \cite{alie1} we checked this conjecture  for $G$ any elliptic curve, in which case  $G^0$ can be taken to be  the complement of the zero section.  
    (This is indeed an easy result; a more subtle result mod $p^2$ for elliptic curves was also proved in \cite{alie1}.)
    Also it is trivial to check that the multiplicative group $G={\mathbb G}_m=GL_1$,
 and hence any split torus, has a  Frobenius lift that is Lie invariant mod $p^m$ for any $m$ (in which case $G^0$ can be taken to be $G$ itself).  The aim of this note is to show that tori are ``essentially" the only linear algebraic groups admitting Frobenius lifts that are Lie invariant mod $p$. To explain this let us state a global consequence (Theorem \ref{immigrationintro} below) of a more refined, local, result (Theorem \ref{immigration}). 

Start with a number field $K$, let $K^a$ be its algebraic closure, let $\cO_K$ be the ring of integers of $K$, and let $G_{\cO_K}$ be an affine group scheme of finite type over $\cO_K$. We denote by $G_{K^a}$ the algebraic group over $K^a$  obtained from $G_{\cO_K}$ via base change. Furthermore, for any unramified prime ${\mathfrak p}$ in $\cO_K$ denote by
$R_{\mathfrak p}=\widehat{\cO_{\mathfrak p}^{\text{ur}}}$
the $p$-adic completion of the maximum unramified extension 
$\cO_{\mathfrak p}^\text{ur}$
of 
the ${\mathfrak p}$-adic completion $\cO_{\mathfrak p}$ of $\cO_K$. Finally denote by $G_{R_{\mathfrak p}}$ the group scheme over $R_{\mathfrak p}$ obtained from $G_{\cO_K}$ via base change. 
So, in particular, if $G_{K^a}$ is a torus then for all but finitely many primes ${\mathfrak p}$ in $\cO_K$ there is a Frobenius lift on  $\widehat{G_{\mathfrak p}}$ that is Lie invariant mod $p^m$ for all $m$. On the other hand we will prove:

 \begin{thm}
 \label{immigrationintro}
 Assume $G_{K^a}$ is connected and is not a torus. Then for all except finitely many primes ${\mathfrak p}$ in $\cO_K$ the group scheme $G:=G_{\mathfrak p}$ over $R:=R_{\mathfrak p}$ satisfies the following property. 
 There is no open subset $G^0\subset G$ 
 with $\widehat{G^0}\neq\emptyset$  such that $\widehat{G^0}$ admits a Frobenius lift  that is Lie invariant mod $p$. \end{thm}
 
 \begin{remark}
 Note that if $K={\mathbb Q}$ and $G$ is any of the  classical split groups 
 $$GL_n,SL_n,Sp_{2n},SO_{2n},SO_{2n+1}$$
 over $\bZ=\cO_K$,
where $n\geq 2$, then the above Theorem will hold for all primes.
\end{remark}

The moral of our Theorem \ref{immigrationintro}, and of results in  \cite{foundations}, Chapter 4, Section 4.11, is that, in the context of our theory in \cite{book, foundations},  there is no ``na\"{i}ve"  arithmetic analogue of translation invariance of vector fields on linear algebraic groups; this state of affairs can be viewed as an a posteriori justification for considering  more ``sophisticated" arithmetic analogues of translation invariance   as in \cite{foundations}. For instance an arithmetic analogue of Lie algebras of algebraic groups was introduced in \cite{foundations, adel2} which does not rely on translation  invariance. Also it is easy to see, cf. \cite{foundations, adel2},  that, in classical differential geometry, translation invariance of vector fields on frame bundles is a {\it consequence} of the conditions defining Chern and Levi-Civit\`{a} connections; on the other hand 
the arithmetic analogues of these conditions 
make sense and therefore  can be taken as a substitute, in arithmetic, for a strengthening of the concept of translation invariance.

\medskip

 The paper is organized as follows: in section 2 we  introduce the relevant concepts and notation; in section 3 we state and prove our main result.

\medskip

{\bf Acknowledgments}. The present work was partially supported by the  IHES and  the Simons Foundation (award 311773). 

\section{Main concepts}

\subsection{Frobenius lifts} Following \cite{alie1} we quickly review some basic notation and terminology pertaining to Frobenius lifts. Let  $p$ be a rational  prime. 
 If $S$ is a ring a {\it Frobenius lift} 
on $S$  is a ring endomorphism $\phi:S\ra S$ whose reduction mod $p$ is the $p$-power 
Frobenius on $S/pS$.

Similarly if $X$ is a scheme or a $p$-adic formal scheme a 
{\it Frobenius lift} on $X$  is an endomorphism $\phi:X\ra X$ whose reduction mod $p$ 
is the $p$-power Frobenius on the reduction of $X$ mod $p$. 
If $X=Spec\ S$ or $X=Spf\ S$ we usually denote by same letter $\phi$ the Frobenius 
lifts on $S$ and $X$. 

Throughout the paper we let $R$ denote any complete discrete valuation ring with maximal ideal generated by $p$ and algebraically closed
residue field $\overline{R}:=R/pR$; 
such an $R$ is uniquely determined up to isomorphism by its residue field (it is the Witt ring on that field) and possesses a unique Frobenius 
lift $\phi:R\ra R$.
For any $R$-algebra $S$ and any scheme or $p$-adic formal scheme $X$ over $R$, Frobenius lifts
on $S$ or $X$ will be assumed compatible with the Frobenius lift on $R$. For any ring $S$ and  scheme $X$  we denote by $\widehat{S}$ and $\widehat{X}$ 
the $p$-adic completions of $S$ and $X$ respectively. 
We also define the K\"{a}hler differentials on the affine formal scheme $\widehat{X}$ over $R$ by
$$\Omega_{\widehat{X}/R}:=\lim_{\leftarrow} \Omega_{X_n/R_n}$$
where $R_n=R/p^nR$, $X_n=X\otimes R_n$, and $\Omega_{X_n/R_n}$ are the usual K\"{a}hler differentials. If $X$ is smooth over $R$ and $\phi$ is a Frobenius lift on $\widehat{X}$  then $\phi$  induces an additive map
\begin{equation}
\label{kedd}
\frac{\phi^*}{p}:\Omega_{\widehat{X}/R}\ra \Omega_{\widehat{X}/R},\ \ \ \omega\mapsto \frac{\phi^*\omega}{p}\end{equation} which,
following \cite{BYM, foundations}, we  view  as an analogue of the Lie derivative on forms in classical differential geometry. 
By the way, the operator \ref{kedd} is a lift to characteristic zero of the Cartier operator and plays a key role in the $p$-adic theory of differential forms. 
 
 \subsection{Lie invariance}
For a  smooth group scheme $G$ over $R$ of relative dimension $d$ let $\omega$ be a column vector of size $D\geq d$ whose entries generate the $R$-module of  left invariant $1$-forms on $G$. Moreover let $G^0\subset G$ be a Zariski open set with non-empty reduction mod $p$. 

 \begin{definition}\label{saul2} A Frobenius lift $\phi$ on $\widehat{G^0}$ is  {\it Lie invariant mod $p^m$}
if there exists an $D\times D$ matrix $\lambda$ with coefficients in $R$ such that:
\begin{equation}
\label{oprah}
\frac{\phi^*}{p}\omega\equiv \lambda \omega\ \ \ \text{mod}\ \ \ p^m.\end{equation}
\end{definition}

The above condition  is independent of the choice of $\omega$; if in addition $\omega$ is a basis for the module of left invariant forms (i.e., $D=d$) the matrix $\lambda$ is, of course, uniquely determined by $\omega$. The case when $G$ is an elliptic curve was analyzed in \cite{alie1} for $m=1,2$, in which case $\lambda$ was related to the Hasse polynomial and some $p$-adic deformations of it. 

   \begin{remark}
Assume, in this remark only, that 
$G$ is a linear algebraic group 
 over the complex numbers. Let
 $\d$ be a ${\mathbb C}$-derivation on the ring $\cO(G)$ (i.e., a {\it vector field} on $GL$)
and let $\text{Lie}_{\d}$ be the induced ${\mathbb C}$-linear endomorphism of $\Omega_{G/{\mathbb C}}$  (the
{\it Lie derivative} along $\d$). Let $\omega$
be a basis for the {\it left invariant} $1$-forms on $G$, i.e., $1$-forms invariant under the automorphisms of $\cO(G)$ given by left translations. Also let us say that $\d$ is {\it right invariant} if it commutes with all automorphisms of $\cO(G)$ given by right translations.
Then the following are equivalent:

1) $\text{Lie}_{\d} \omega=0$;

2) $\d$ is right invariant.

\medskip

We would like to see condition \ref{oprah}  as an arithmetic analogue
of the Condition 1); cf. also \cite{BYM, euler, alie1}. Note that, in \cite{foundations}, Chapter 4, Section 4.11, we proved that, for $G=GL_n$ over $R$,  there are no Frobenius lifts satisfying certain arithmetic analogues of Condition 2). The present paper can be viewed as being devoted to proving that there are no Frobenius lifts satisfying an appropriate arithmetic analogue of condition 1). Taken together these results can be viewed as a justification for the introduction
of a number of substitutes for invariance of vector fields in \cite{foundations}. 
 \end{remark}

 \subsection{Linear algebraic groups}
For our next discussion let us consider
 the multiplicative, additive, and general linear  group schemes over an arbitrary ring $A$:
$$
\begin{array}{rcl}
{\mathbb G}_m & = & Spec\ A[s,s^{-1}],\\
\ & \ & \ \\
{\mathbb G}_a & = & Spec\ A[t],\\
\ & \ & \ \\
GL_N & = & Spec\ A[x,\det(x)^{-1}],\ \ x=(x_{ij})_{i,j=1}^N.
\end{array}
$$
As usual we view $GL_N$ as embedded in 
$${\mathfrak g}{\mathfrak l}_N:=Spec\ A[x],\ \ x=(x_{ij})_{i,j=1}^N.$$
Recall that a matrix $u\in GL_N(A)$ is called {\it unipotent} if $u-1$
is nilpotent in ${\mathfrak g}{\mathfrak l}_N(A)$, where $1=1_N$ is the identity matrix.
Also recall that if $G$ is an affine group scheme of finite type over $A$ then an element $u\in G(A)$ is called {\it unipotent} if there exists an embedding $G\subset GL_N$
such that $u-1$ is nilpotent in ${\mathfrak g}{\mathfrak l}_N(A)$

\begin{definition}
Let $G$ be an affine group scheme of finite type and let   $u\in G(A)$ be a unipotent element.
The {\it index} of $u$ is the minimum of all integers $n\geq 1$ having the property that there exists an embedding $G\subset GL_N$
such that $(u-1)^n=0$ in ${\mathfrak g}{\mathfrak l}_N(A)$.\end{definition}

\begin{remark}\label{vinno}\

1) If $G$ is a  torus over an integral domain $A$ (i.e., $G$ is isomorphic over the algebraic closure of the fraction field of $A$ to a power of ${\mathbb G}_m$) 
 then $G(A)$ has no non-trivial unipotent elements.

2) Let $G$ over $A=\bZ$ be either a power of ${\mathbb G}_a$ or any of the classical groups   
$$GL_n,  SL_n, Sp_{2n}, SO_{2n}, SO_{2n+1}.$$
Here the $Sp_N$ and $SO_N$ groups are the split forms of the symplectic and orthogonal groups, i.e., the   identity components of the subgroup schemes of the corresponding $GL_N$
 defined by the entries of $x^tqx-q$, where
 $$q=\left(\begin{array}{cc}0 & 1_n\\-1_n & 0\end{array}\right),\ \ \left(\begin{array}{ll}0 & 1_n\\1_n & 0\end{array}\right),\ \ \left(\begin{array}{lll}1 & 0 & 0\\
 0& 0 & 1_n\\0 & 1_n & 0\end{array}\right),$$
 respectively. Then it is well known that $G(\bZ)$ contains  unipotent elements of index $2$.  Cf. \cite{foundations}, pp. 54 and 180. 
 
 3) Recall the following situation considered in the Introduction. Let $K$ be  a number field, let $K^a$ be its algebraic closure, let $\cO_K$ be the ring of integers of $K$, and let $G_{\cO_K}$ be an affine group scheme of finite type over $\cO_K$. We denote by $G_{K^a}$ the algebraic group over $K^a$  obtained from $G_{\cO_K}$ via base change. Furthermore, for any unramified prime ${\mathfrak p}$ in $\cO_K$ denote by
$R_{\mathfrak p}=\widehat{\cO_{\mathfrak p}^{\text{ur}}}$
the $p$-adic completion of the maximum unramified extension 
$\cO_{\mathfrak p}^\text{ur}$
of 
the ${\mathfrak p}$-adic completion $\cO_{\mathfrak p}$ of $\cO_K$. Finally denote by $G_{R_{\mathfrak p}}$ the group scheme over $R_{\mathfrak p}$ obtained from $G_{\cO_K}$ via base change.  Assume now  $G_{K^a}$  is connected but not a torus.
 Then, by \cite{hum}, p. 161, $G_{K^a}(K^a)$ contains a non-trivial unipotent element. As a consequence,
 for all but finitely many primes ${\mathfrak p}$ in $\cO_K$ the group $G=G_{R_{\mathfrak p}}$ over $R:=R_{\mathfrak p}$ has the property that  $G(R)$ contains a non-trivial  unipotent element of index $\leq p$.
\end{remark}

\section{Main result}
The moral of our paper will be that, for affine group schemes,  Frobenius lifts that are Lie invariant mod $p$ exist only in very special cases. First note:

\begin{remark}
 Let $G$ be a split torus over $R$, i.e., a power of ${\mathbb G}_m=Spec\ R[s,s^{-1}]$.  The Frobenius lift $\phi$ on $G$ induced  by $\phi(s)=s^p$ 
 has, of course, the property that
 $\frac{\phi^*}{p}\omega=\omega$ where $\omega$ has entries $\frac{ds}{s}$. In particular
 the  Frobenius lift on $\widehat{G}$ induced by $\phi$
 is Lie invariant mod $p^m$ for all $m$.  \end{remark}
 
Our main result is the following:

\begin{thm}
\label{immigration}
 Let $G$ be an affine  smooth group scheme over $R$. Assume
 there exists an open subset $G^0\subset G$ 
 with $\widehat{G^0}\neq\emptyset$  such that $\widehat{G^0}$ admits a Frobenius lift  Lie invariant mod $p$. Then $G(R)$ contains no non-trivial unipotent elements of index $\leq p$.
\end{thm}

Theorem \ref{immigration} and  Remark \ref{vinno} imply Theorem \ref{immigrationintro} in the Introduction; they also imply:
    
\begin{corollary}\label{gugu}
Let $G$ over $R$ be either a power of ${\mathbb G}_a$ or one of the split classical groups 
$$GL_n,SL_n,Sp_{2n},SO_{2n},SO_{2n+1}$$ 
where $n\geq 2$. Then
there is  no open set $G^0\subset G$ with
    $\widehat{G^0}\neq\emptyset$  such that $\widehat{G^0}$ admits a Frobenius lift that is Lie invariant mod $p$.
\end{corollary}

The rest of this note is devoted to the proof of Theorem 
   \ref{immigration}.  We need the following easy:

\begin{lemma}
\label{denada}
Let $G$ be a smooth affine group scheme over $R$. 
Assume $G(R)$ contains
a non-trivial unipotent element $u$ of index $\leq p$.  Then there exists a group homomorphism  $\iota:{\mathbb G}_a\ra G$ over $R$ and an embedding $G\subset GL_N$ over $R$ such that 
the composition 
$j:{\mathbb G}_a\stackrel{\iota}{\longrightarrow} G\subset GL_N$
is given by
 \begin{equation}\label{suc}
 j^*:R[x,\det(x)^{-1}]\ra R[t],\ \ j^*x\equiv \left(\begin{array}{cc} 1_2+t\nu_2 & \star\\
   0 & \star\end{array}\right), \ \ \ \text{mod}\ \ p\end{equation}
   where $1_2,\nu_2$ are the $2\times 2$ matrices,
   $$1_2:=\left(\begin{array}{cc}
   1 & 0 \\ 0 & 1\end{array}\right),\ \ \nu_2:=\left(\begin{array}{cc}
   0 & 1 \\ 0 & 0\end{array}\right).$$ \end{lemma}
   
   {\it Proof}. 
   Let $G\subset GL_N$ be an embedding such that  $(u-1)^p=0$. 
   Write $u=1+p^e\nu$ where $e\geq 0$ and $\nu\not\equiv 0$ mod $p$.
   Since $\nu^p=0$ we have that
    $\beta:=p^{-e}\text{log}(u)$ has coefficients in $R$,
   $\beta^p=0$, 
    and hence  $\text{exp}(t\beta)$ has coefficients in $R[t]$. So we may define
   a homomorphism 
   $$\iota:{\mathbb G}=Spec\ R[t]\ra GL_N=Spec\ R[x,\det(x)^{-1}],\ \ x\mapsto \text{exp}(t\beta).$$
   We claim that this homomorphism factors through a homomorphism $\iota:{\mathbb G}_a \ra G$. 
   Indeed, let $F_j$ be the polynomials (with coefficients in $R$) 
 defining the ideal of $G$ in $GL_n$. 
 We know that $F_j$ vanish on $u$ and we want to check they vanish
 on $\text{exp}(t\beta)$. 
 Let $f_j(t):=F_j(\text{exp}(t\beta))\in R[t]$. Since
 $\text{exp}(p^ek\beta)=u^k$ it follows that
  that $f_j(p^ek)=0$ for all
 $k\in \bZ$. So $f_j(t)=0$ and our claim is proved.
   
   Let us denote by 
   an overline the reduction mod $p$; in particular
   $\overline{\nu}\in {\mathfrak g}{\mathfrak l}_N(\overline{R})$ denotes the image of $\nu$.
   By the theory of Jordan normal form, there is a matrix $\overline{v}\in GL_N(\overline{R})$ such that $\overline{v}\cdot \overline{\nu}\cdot \overline{v}^{-1}=(\overline{\nu'}_{ij})$,  $\nu'_{i,i+1}\in \{0,1\}$, $\nu'_{12}=1$, and all other $\nu'_{ij}=0$. Lift $\overline{v}$ to a matrix $v\in GL_N(R)$. Changing the embedding $G\subset GL_N$ via conjugation by $v$ we may assume $\nu\equiv \nu'$ mod $p$. Set $u':=1+p^e\nu'$ and $\beta':=p^{-e}\text{log}(u')$.
  Since $\nu^p=(\nu')^p=0$ it follows that
  $$\beta\equiv \beta'\ \ \ \text{mod}\ \ p.$$
  Since $\beta^p=(\beta')^p=0$ it follows that
  $$\text{exp}(t\beta)\equiv \text{exp}(t\beta')\ \ \ \text{mod}\ \ p.$$
  On the other hand
  $$
 \text{exp}(t\beta')\equiv  \left(\begin{array}{cc} 1_2+t\nu_2 & \star\\
   0 & \star\end{array}\right), \ \ \ \text{mod}\ \ p.
  $$
   Then \ref{suc} holds and we are done.   
    \qed

    \bigskip   
   
   {\it Proof of Theorem \ref{immigration}}.    
   Assume the hypotheses of the Theorem and assume $G(R)$ contains a non-trivial unipotent element of index $\leq p$; we will derive a contradiction. By Lemma \ref{denada} there exists a homomorphism
     $\iota: {\mathbb G}_a\ra  G$  and a closed embedding $G\subset GL_N$ 
   such that the composition ${\mathbb G}_a\ra GL_N$ is given by \ref{suc}.
So if we set $\xi=(\xi_{ij})$, where  $\xi_{ij}\in \cO(G)$ are the images of $x_{ij}$, then
   $$\iota^*\xi\equiv  \left(\begin{array}{cc} 1_2+t\nu_2 & \star\\
   0 & \star\end{array}\right) \ \ \ \text{mod}\ \ \ p.$$
  
  Now for any $g\in G(R)$ denote by $L_g:G\ra G$, $L_g:\widehat{G}\ra \widehat{G}$ the left translation by $g$. Let $\phi$ be a Frobenius lift on $\widehat{G^0}$ that is Lie invariant mod $p$.
  Clearly $L_g\circ \phi \circ L_{g^{-1}}$ is a Frobenius lift
  on $L_g(\widehat{G^0})$ which is still Lie invariant mod $p$; hence replacing $\phi$
  by $L_g\circ \phi \circ L_{g^{-1}}$ 
  and $G^0$ by $L_g(G^0)$, for an appropriate $g$, 
  we may assume ${\mathbb G}_a^0:=\iota^{-1}(\widehat{G^0})$ is non-empty.  We continue to denote by $\iota:{\mathbb G}_a^0\ra \widehat{G^0}$ the restriction of $\iota$.
   
  For any $N\times N$ matrix $M$ we denote by $M^{\text{col}}$ the column of size $N^2$
  obtained by concatenating the columns of $M$. 
 Set $$\omega:=(\xi^{-1}d\xi)^{\text{col}}.$$ 
   The entries of $\omega$ generate  the $R$-module of left invariant $1$-forms on $G$. 
  Write
  $$\phi(\xi)=\xi^{(p)}+pZ$$
  where $\xi^{(p)}=(\xi_{ij}^p)$ and $Z=(Z_{ij})$ is an $N\times N$ matrix with entries in $\cO(\widehat{G^0})$. By our assumptions we have
  \begin{equation}
  \label{ionescuuuu}
  \frac{\phi^*}{p}\omega\equiv \lambda\cdot \omega\ \ \ \text{mod}\ \ \ p,
  \end{equation}
 for some $N^2\times N^2$ matrix $\lambda$ with coefficients in $R$.
  The congruence \ref{ionescuuuu} yields:
  \begin{equation}
  \label{jean}
  ((\xi^{(p)})^{-1}(\eta+dZ))^{\text{col}}\equiv \lambda \cdot (\xi^{-1}d\xi)^{\text{col}}\ \ \ \text{mod}\ \ p,\end{equation}
   where
   $$\eta=(\eta_{ij}),\ \ \eta_{ij}=\xi_{ij}^{p-1}d\xi_{ij}.$$
  Set $z_{ij}=\iota^*Z_{ij}$, $z=(z_{ij})$.    Then we have
   the following congruences in $\cO({\mathbb G}_a^0)dt$:
   $$
   \begin{array}{rcll}
   \iota^*((x^{(p)})^{-1}) & \equiv & \left(\begin{array}{cc}
   1_2-t^p\nu_2 & \star\\
   0 & \star\end{array}\right) & \text{mod}\ \ p\\
   \ & \ & \ & \ \\
   \iota^*\eta & \equiv & \left(\begin{array}{cc} t^{p-1}dt\cdot \nu_2 & \star\\
   0 & \star\end{array}\right) & \text{mod}\ \ p\\
   \ & \ & \ & \ \\
   \iota^*(dZ) &\equiv & dz & \text{mod}\ \ p.
   \end{array}
   $$   
   Finally $\iota^* (\xi^{-1}d\xi)$ mod $p$ has left invariant entries (because $\iota$ is a homomorphism) hence these entries are congruent mod $p$ to $R$-multiples of $dt$.
   Taking $\iota^*$ in \ref{jean}, picking out the $(N+1)$-th entry, and dividing by $dt$, we get the following congruence in $\cO({\mathbb G}_a^0)$, and hence in $\widehat{R((t))}$:  
    \begin{equation}\label{rosa}
   t^{p-1}+\frac{dz_{12}}{dt}-t^p\frac{dz_{22}}{dt}\equiv \mu\ \ \ \text{mod}\ \ p,\end{equation}
   for some $\mu\in R$.
   Viewing $z_{12},z_{22}$ as series in $t$ we see that
   the coefficients of $t^{p-1}$ in the left and right hand sides of \ref{rosa} are congruent to $1$ and $0$ mod $p$ respectively, a contradiction; this ends our proof.
     \qed

\end{document}